\documentclass[reqno, 11pt, oneside]{amsart}

\usepackage{latexsym, amsmath,amssymb}
\usepackage{hyperref}
\hypersetup{colorlinks=true, citecolor=cyan, linkcolor=blue}

\usepackage[initials, alphabetic,msc-links, nobysame]{amsrefs}
\usepackage{graphicx}
\usepackage{xcolor}
\usepackage{enumitem}

\usepackage[utf8]{inputenc}

\usepackage{amsmath}

\def\XXint#1#2#3{{\setbox0=\hbox{$#1{#2#3}{%
\int}$ }
\vcenter{\hbox{$#2#3$ }}\kern-.6\wd0}}

\renewcommand{\epsilon}{\varepsilon}
\newtheorem{theorem}{Theorem}

\newtheorem{lemma}[theorem]{Lemma}
\newtheorem{corollary}[theorem]{Corollary}

\newtheorem{proposition}[theorem]{Proposition}
\newtheorem{definition}[theorem]{Definition}
\newtheorem{remark}[theorem]{Remark}
\newcommand{\bth}{\begin{theorem}}
\newcommand{\ble}{\begin{lemma}}
\newcommand{\bcor}{\begin{corr}}

\newcommand{\bdeff}{\begin{deff}}

\newcommand{\bprop}{\begin{proposition}}
\newcommand{\ele}{\end{lemma}}
\newcommand{\ecor}{\end{corr}}
\newcommand{\edeff}{\end{deff}}

\newcommand{\eprop}{\end{proposition}}

\newcommand{\supp}{\text{supp }}
\renewcommand{\Pi}{\varPi}

\renewcommand{\epsilon}{\varepsilon}

\newcommand{\R}{{\mathbb R}}

\usepackage{esint}

\usepackage{comment}

\allowdisplaybreaks

\def\vint_#1{\mathchoice%
        {\mathop{\kern 0.2em\vrule width 0.6em height 0.69678ex depth -0.58065ex
                \kern -0.8em \intop}\nolimits_{\kern -0.4em#1}}%
        {\mathop{\kern 0.1em\vrule width 0.5em height 0.69678ex depth -0.60387ex
                \kern -0.6em \intop}\nolimits_{#1}}%
        {\mathop{\kern 0.1em\vrule width 0.5em height 0.69678ex
            depth -0.60387ex
                \kern -0.6em \intop}\nolimits_{#1}}%
        {\mathop{\kern 0.1em\vrule width 0.5em height 0.69678ex depth -0.60387ex
                \kern -0.6em \intop}\nolimits_{#1}}}
\def\vintslides_#1{\mathchoice%
        {\mathop{\kern 0.1em\vrule width 0.5em height 0.697ex depth -0.581ex
                \kern -0.6em \intop}\nolimits_{\kern -0.4em#1}}%
        {\mathop{\kern 0.1em\vrule width 0.3em height 0.697ex depth -0.604ex
                \kern -0.4em \intop}\nolimits_{#1}}%
        {\mathop{\kern 0.1em\vrule width 0.3em height 0.697ex depth -0.604ex
                \kern -0.4em \intop}\nolimits_{#1}}%
        {\mathop{\kern 0.1em\vrule width 0.3em height 0.697ex depth -0.604ex
                \kern -0.4em \intop}\nolimits_{#1}}}

\newcommand{\aveint}[2]{\mathchoice%
        {\mathop{\kern 0.2em\vrule width 0.6em height 0.69678ex depth -0.58065ex
                \kern -0.8em \intop}\nolimits_{\kern -0.45em#1}^{#2}}%
        {\mathop{\kern 0.1em\vrule width 0.5em height 0.69678ex depth -0.60387ex
                \kern -0.6em \intop}\nolimits_{#1}^{#2}}%
        {\mathop{\kern 0.1em\vrule width 0.5em height 0.69678ex depth -0.60387ex
                \kern -0.6em \intop}\nolimits_{#1}^{#2}}%
        {\mathop{\kern 0.1em\vrule width 0.5em height 0.69678ex depth -0.60387ex
                \kern -0.6em \intop}\nolimits_{#1}^{#2}}}

\newcommand{\loc}{\mathrm{loc}}

\def\XXint#1#2#3{{\setbox0=\hbox{$#1{#2#3}{\int}$}
    \vcenter{\hbox{$#2#3$}}\kern-.5\wd0}}

\newcommand{\vertiii}[1]{{\left\vert\kern-0.25ex\left\vert\kern-0.25ex\left\vert #1 
    \right\vert\kern-0.25ex\right\vert\kern-0.25ex\right\vert}}

\newcommand{\vertii}[1]{{\left\vert\kern-0.25ex\left\vert\kern-0.25ex  #1 
    \kern-0.25ex\right\vert\kern-0.25ex\right\vert}}

\makeatletter
\@namedef{subjclassname@2020}{\textup{2020} Mathematics Subject Classification}
\makeatother

\usepackage{mathabx}
\usepackage{geometry}

\begin{document}

\title{Logarithmic Poincaré--Sobolev inequalities}

\author[A. Claros]{Alejandro Claros}

\address[A. Claros]{BCAM -- Basque Center for Applied Mathematics, Bilbao, Spain\newline
Universidad del Pa\'is Vasco / Euskal Herriko Unibertsitatea (UPV/EHU), Bilbao, Spain}

\email{aclaros@bcamath.org, aclaros003@ikasle.ehu.eus}

\author[I. Gardeazabal-Guti\'errez]{Iker Gardeazabal-Guti\'errez}

\address[I. Gardeazabal-Guti\'errez]{Universidad del Pa\'is Vasco / Euskal Herriko Unibertsitatea \newline (UPV/EHU), Bilbao, Spain}

\email{iker.gardeazabal@ehu.eus}

\author[C. P\'erez]{Carlos P\'erez}

\address[C. P\'erez]{BCAM -- Basque Center for Applied Mathematics, Bilbao, Spain\newline
Universidad del Pa\'is Vasco / Euskal Herriko Unibertsitatea (UPV/EHU), Bilbao, Spain\newline
Ikerbasque, Bilbao, Spain}

\email{cperez@bcamath.org}

\thanks{A. Claros is supported by the Basque Government through the BERC 2022-2025 program, by the Ministry of Science and Innovation through Grant PRE2021-099091 funded by BCAM Severo Ochoa accreditation CEX2021-001142-S/MICIU/AEI/10.13039/501100011033 and by ESF+, and by the project PID2023-146646NB-I00 funded by MICIU/AEI/10.13039/501100011033 and by ESF+. 
\\ \hspace*{1.5em} I. Gardeazabal-Gutiérrez is supported by EHU project PES 22/50. 
\\ \hspace*{1.5em}C. P\'erez is supported by the Spanish government through the grant PID2023-146646NB-I00 and by Severo Ochoa accreditation CEX2021-001142-S, both at BCAM, and also by the Basque Government through grant  IT1875-26 at the University of the Basque Country and by the BERC programme 2022-2025 at BCAM}

\subjclass[2020]{Primary 46E35, 26D10; Secondary 26A33, 42B35}

\keywords{Logarithmic Sobolev inequalities, Poincaré--Sobolev inequalities, fractional Sobolev spaces, Bourgain--Brezis--Mironescu formula, self-improving inequalities, Nash inequalities}

\begin{abstract}
We establish local logarithmic Poincaré--Sobolev inequalities on cubes in $\R^n$ in both the classical and fractional settings. These estimates follow by combining Poincaré--Sobolev inequalities with an elementary consequence of Jensen's inequality on probability spaces. The same mechanism yields logarithmic bounds from the self-improving properties of generalized Poincaré inequalities. For $p>1$, we recover the corresponding global logarithmic Sobolev inequalities on $\R^n$ by applying the local estimates on expanding cubes. In the fractional setting, our inequalities incorporate the Bourgain--Brezis--Mironescu factor $1-\delta$ and recover the first-order counterparts   as $\delta\to1^-$. As a further consequence, we derive both classical and fractional local $L^p$ Nash inequalities.

\end{abstract}

\maketitle

\section{Introduction and statements of the main results}

The classical Euclidean logarithmic Sobolev inequality states that, for $n\geq1$ and $f\in C_c^\infty(\R^n)$ such that
\begin{equation*}
	\int_{\R^n}|f(x)|^2\,dx=1,
\end{equation*}
we have
\begin{equation}\label{eq:log-global-L2}
	\int_{\R^n}|f(x)|^2\log|f(x)|^2\,dx\le\frac{n}{2}\log \left(\frac{2}{\pi e n}\int_{\R^n}|\nabla f(x)|^2\,dx\right).
\end{equation}
This inequality is equivalent to the Gaussian logarithmic Sobolev inequality of Gross \cite{Gross} (see Weissler \cite{Weissler}), and the constant $\frac{2}{\pi e n}$ is optimal. The equality cases in the logarithmic Sobolev inequalities for Lebesgue and Gaussian measures are characterized by Carlen \cite{Carlen}. The Euclidean $L^2$ version also arises as a limiting case of a sharp family of Gagliardo--Nirenberg inequalities obtained by Del Pino and Dolbeault~\cite{DPD02}. There is also a classical $L^p$ version of \eqref{eq:log-global-L2}: for $1<p<\infty$ and $f\in C_c^\infty(\R^n)$ with
\begin{equation*}
	\int_{\R^n}|f(x)|^p\,dx=1,
\end{equation*}
we have
\begin{equation}\label{eq:log-global-Lp}
	\int_{\R^n}|f(x)|^p\log|f(x)|^p\,dx \le\frac np\log\left(\alpha_{n,p}\int_{\R^n}|\nabla f(x)|^p\,dx\right),
\end{equation}
where $\alpha_{n,p}$ denotes the optimal constant. Its value and explicit extremizers were obtained by Del Pino and Dolbeault \cite{DPD03} for $1<p<n$ and by Gentil \cite{Gentil} for all $p>1$.

\begin{remark}\label{rem:scaling}
	The coefficient $\frac np$ in \eqref{eq:log-global-Lp} is given by scaling. For $\lambda>0$, the function $f_\lambda(x)=\lambda^{n/p}f(\lambda x)$ preserves the $L^p$ norm of $f$, while
	\begin{equation*}
		\int_{\R^n}|\nabla f_\lambda(x)|^p\,dx=\lambda^p\int_{\R^n}|\nabla f(x)|^p\,dx.
	\end{equation*}
Hence, replacing $f$ by $f_\lambda$ increases the left-hand side of~\eqref{eq:log-global-Lp} by $n \log \lambda$, and the right-hand side increases by the exact same amount, since $\frac np \log \lambda^p = n \log \lambda$.
\end{remark}

The purpose of this paper is to obtain local versions of these inequalities on cubes. On the one hand, a local formulation should measure the oscillation of $f$ around its average, since the gradient does not change when a constant is added to $f$. On the other hand, it should take into account the size of the cube and recover the global inequality when the cube increases to the whole space. Recently, Abdo and Lee \cite{AL25} proved logarithmic Sobolev inequalities on bounded smooth domains through interpolation estimates and applied them to drift-diffusion equations.

Given a cube $Q\subset\R^n$, we denote its side length by $\ell(Q)$ and write
\begin{equation*}
	f_Q:=\frac{1}{|Q|}\int_Q f(x)\,dx.
\end{equation*}
Throughout the paper, functions are real-valued and we use the convention $0\log0=0$. Our first main result is the following local logarithmic Poincaré--Sobolev inequality.

\begin{theorem}\label{thm:intro-local}
	Let $1\le p<n$. There exists a constant $C_{n,p}>0$ such that, for every cube $Q\subset\R^n$ and every nonconstant $f\in W^{1,p}(Q)$, we have
	\begin{equation}\label{eq:intro-local}
	\begin{aligned}
		&\frac{1}{|Q|}\int_Q|f(x)-f_Q|^p\log \left(\frac{|f(x)-f_Q|^p}{\frac{1}{|Q|}\int_Q|f(y)-f_Q|^p\,dy}\right)dx\\
		&\quad\le\frac np\left(\frac{1}{|Q|}\int_Q|f(x)-f_Q|^p\,dx\right)\log \left(\frac{C_{n,p}^p\ell(Q)^p\frac{1}{|Q|}\int_Q|\nabla f(x)|^p\,dx}{\frac{1}{|Q|}\int_Q|f(y)-f_Q|^p\,dy}\right).
	\end{aligned}
	\end{equation}
\end{theorem}

The probability measure $\frac{dx}{|Q|}$ and the factor $\ell(Q)^p$ provide the natural local normalization of the inequality. Both sides of \eqref{eq:intro-local} are homogeneous of degree $p$ with respect to multiplication of $f$ by scalars. Here nonconstant means that $f$ is not constant almost everywhere on $Q$, which is equivalent to $\frac{1}{|Q|}\int_Q|f(y)-f_Q|^p\,dy>0$. If $f$ is constant, both sides can be interpreted as zero.

We will show in Proposition~\ref{prop:local-implies-global} that \eqref{eq:intro-local} implies the classical global logarithmic Sobolev inequality for $1<p<n$. More precisely, applying \eqref{eq:intro-local} on an increasing family of cubes and passing to the limit gives \eqref{eq:log-global-Lp} with $C_{n,p}^p$ in place of $\alpha_{n,p}$. The terms involving the side length cancel because $\frac np\log\ell(Q)^p=n\log\ell(Q)$. Thus \eqref{eq:intro-local} combines the local form with the scaling needed to recover the classical global inequality. The approach in \cites{DPD02,DPD03}, based on optimal Gagliardo--Nirenberg inequalities on $\R^n$ and their connection with nonlinear diffusion equations, depends on the structure of the whole space and does not seem to adapt to cubes. Our argument only uses the Poincaré--Sobolev inequality and also gives a fractional version with the Bourgain--Brezis--Mironescu factor $(1-\delta)$ (Theorem~\ref{thm:intro-fractional} below).

\begin{remark}
	The constant $C_{n,p}$ can be taken to be any constant for which the classical Poincaré--Sobolev inequality
	\begin{equation*}
		\left(\frac{1}{|Q|}\int_Q|f(x)-f_Q|^{p^*}\,dx\right)^{\frac{1}{p^*}}\le C_{n,p}\ell(Q)\left(\frac{1}{|Q|}\int_Q|\nabla f(x)|^p\,dx\right)^{\frac{1}{p}}
	\end{equation*}
	holds, where $p^*=\frac{np}{n-p}$. The coefficient $\frac np$ follows from the identity $\frac{p^*}{p^*-p}=\frac np$. Our argument does not recover the optimal global constant $\alpha_{n,p}$; it gives the admissible constant $C_{n,p}^p$, and hence $\alpha_{n,p}\le C_{n,p}^p$ for $1<p<n$.
\end{remark}

The proof of Theorem \ref{thm:intro-local} uses a simple estimate on probability spaces. By Jensen's inequality, the $L^p$ norm of a nonnegative function controls an $L\log L$-type integral of its $q$-th power for $q<p$. Applying this estimate to $|f-f_Q|$, with the normalized measure on $Q$, gives the logarithmic inequality from the classical Poincaré--Sobolev estimate. The same argument applies whenever a suitable higher integrability estimate for the oscillation is available. In particular, it can be combined with the self-improvement of generalized Poincaré inequalities developed in \cites{FPW98,PR,LLO,ClarosPerezZheng} (see Theorem~\ref{thm:log-selfimproving}).

Our second main result is a fractional version of Theorem \ref{thm:intro-local}. For $0<\delta<1$ and $1\le p<\frac{n}{\delta}$, we set
\begin{equation*}
	p_\delta^*:=\frac{np}{n-\delta p}.
\end{equation*}
We use the fractional Sobolev space $W^{\delta,p}(Q)$ of functions $f\in L^p(Q)$ such that
\begin{equation*}
	\int_Q\int_Q\frac{|f(x)-f(y)|^p}{|x-y|^{n+\delta p}}\,dy\,dx<\infty.
\end{equation*}
For background on fractional Sobolev spaces, we refer to \cite{Hitchhiker}. Our proof relies on the fractional Poincaré--Sobolev inequality with the sharp factor $(1-\delta)^{1/p}$ due to Bourgain, Brezis, and Mironescu~\cite{BBM1}. 

\begin{theorem}\label{thm:intro-fractional}
	Let $0<\delta<1$ and $1\le p<\frac{n}{\delta}$. There exists a constant $C_{n,p,\delta}>0$ such that, for every cube $Q\subset\R^n$ and every nonconstant $f\in W^{\delta,p}(Q)$, we have
	\begin{equation}\label{eq:intro-fractional}
	\begin{aligned}
		&\frac{1}{|Q|}\int_Q|f-f_Q|^p\log \left(\frac{|f-f_Q|^p}{\frac{1}{|Q|}\int_Q|f-f_Q|^p\,dy}\right)dx\\
		&\quad\le\frac{n}{\delta p}\left(\frac{1}{|Q|}\int_Q|f-f_Q|^p\,dx\right) \log \left(\frac{C_{n,p,\delta}^p(1-\delta)\ell(Q)^{\delta p}\frac{1}{|Q|}\int_Q\int_Q\frac{|f(x)-f(y)|^p}{|x-y|^{n+\delta p}}\,dy\,dx}{\frac{1}{|Q|}\int_Q|f-f_Q|^p\,dy}\right).
	\end{aligned}
	\end{equation}
\end{theorem}

The coefficient $\frac{n}{\delta p}$ is obtained from $\frac{p_\delta^*}{p_\delta^*-p}=\frac{n}{\delta p}$. As in the classical case, the local inequality implies a global one. We will prove in Proposition~\ref{prop:fractional-local-implies-global} that, for $1<p<\frac{n}{\delta}$ and $f\in C_c^\infty(\R^n)$ satisfying $\int_{\R^n}|f(x)|^p\,dx=1$, we have
\begin{equation}\label{eq:intro-global-fractional}
	\int_{\R^n}|f(x)|^p\log|f(x)|^p\,dx\le\frac{n}{\delta p}\log \left(C_{n,p,\delta}^p(1-\delta)\int_{\R^n}\int_{\R^n}\frac{|f(x)-f(y)|^p}{|x-y|^{n+\delta p}}\,dy\,dx\right).
\end{equation}
The passage to \eqref{eq:intro-global-fractional} uses the same cancellation of the side lengths, since $\frac{n}{\delta p}\log\ell(Q)^{\delta p}=n\log\ell(Q)$.

\begin{remark}\label{rem:fractional-scaling}
	The computation in Remark~\ref{rem:scaling} also applies to \eqref{eq:intro-global-fractional}. In this case, the seminorm of $f_\lambda$ scales as $\lambda^{\delta p}$, and the coefficient $\frac{n}{\delta p}$ satisfies $\frac{n}{\delta p}\log\lambda^{\delta p}=n\log\lambda$.
\end{remark}

The fractional result recovers Theorem~\ref{thm:intro-local} as $\delta\to1^-$. For fixed $1\le p<n$, the constants $C_{n,p,\delta}$ can be chosen uniformly bounded for $\frac12\le\delta<1$. Indeed, in the unweighted case the constant in \cite{MPW}*{Theorem~5.7} is $C_n\,p_\delta^*\,\delta^{-(1-\frac1p)}\le 2C_n\frac{np}{n-p}$. Hence we may take $C_{n,p,\delta}\le4C_n\frac{np}{n-p}$.
On the other hand, the Bourgain--Brezis--Mironescu formula \cite{BBM2} gives, for every $f\in W^{1,p}(Q)$,
\begin{equation*}
	\lim_{\delta\to1^-}(1-\delta)\int_Q\int_Q\frac{|f(x)-f(y)|^p}{|x-y|^{n+\delta p}}\,dy\,dx=K_{n,p}\int_Q|\nabla f(x)|^p\,dx,
\end{equation*}
where
\begin{equation*}
	K_{n,p}:=\frac1p\int_{\mathbb S^{n-1}}|\omega\cdot e|^p\,d\sigma(\omega),
\end{equation*}
$e$ is any unit vector and $d\sigma$ is the usual surface measure. The formula is stated in \cite{BBM2} for smooth bounded domains. By Bal, Mohanta, and Roy \cite{KKP}*{Theorems~1 and~2}, it holds on every extension domain (for $1<p<\infty$ and for $p=1$, respectively), and in particular on cubes. Letting $\delta\to1^-$ in \eqref{eq:intro-fractional}, we obtain \eqref{eq:intro-local} with a constant depending only on $n$ and $p$.

Finally, in Section~\ref{sec:nash}, we adapt  an argument of Beckner~\cite{Beckner} to derive classical and fractional local $L^p$ Nash inequalities on cubes from Theorems~\ref{thm:intro-local} and~\ref{thm:intro-fractional} (see Propositions~\ref{prop:local-Nash} and~\ref{prop:fractional-Nash}). These estimates incorporate the Bourgain--Brezis--Mironescu factor $1-\delta$ in the fractional regime and yield a local version of the classical Nash inequality when $p=2$ (Corollary~\ref{cor:local-Nash-L2}).

In a forthcoming work, we will study weighted versions of these inequalities, including their fractional counterparts, with an explicit dependence on the Muckenhoupt constant. Since our argument relies only on a Poincaré--Sobolev inequality, it is natural to expect a similar approach to apply to domains supporting such inequalities, such as John domains.

\section{Preliminaries}\label{sec:preliminaries}

We begin with an elementary estimate on probability spaces which will be the key to our approach. This estimate follows from Jensen's inequality and allows us to pass from higher integrability to logarithmic inequalities. We then recall the self-improving properties of generalized Poincaré inequalities and apply the lemma to obtain logarithmic estimates for the oscillation on cubes.

\begin{lemma}\label{lem:probability}
	Let $(X,\mu)$ be a probability space, let $0<q<p<\infty$, and let $g\geq0$ be such that $g\in L^p(X,\mu)$.  If $g$ does not vanish almost everywhere, then
	\begin{equation*}
		\int_X g^q \log \left(\frac{g^q}{\int_X g^q\,d\mu}\right)d\mu \le\frac{p}{p-q}\left(\int_X g^q\,d\mu\right) \log \left(\frac{\left(\int_X g^p\,d\mu\right)^{q/p}} {\int_X g^q\,d\mu}\right).
	\end{equation*}
	In particular, if $\bigl(\int_X g^p\,d\mu\bigr)^{\frac{1}{p}}\le A$, then
	\begin{equation*}
		\int_X g^q\log \left(\frac{g^q}{\int_X g^q\,d\mu}\right)d\mu \le\frac{p}{p-q}\left(\int_X g^q\,d\mu\right) \log \left(\frac{A^q}{\int_X g^q\,d\mu}\right) \le\frac{p}{e(p-q)}A^q.
	\end{equation*}
\end{lemma}

\begin{proof}
	Since $0<q<p$ and $\mu(X)=1$, we have $0<\int_X g^q\,d\mu<\infty$. We define
	\begin{equation*}
		d\sigma:=\frac{g^q}{\int_X g^q\,d\mu}\,d\mu.
	\end{equation*}
	Hence, we have $\sigma(X)=1$, and the set $\{g>0\}$ has full $\sigma$-measure. Set $\theta=\frac{p}{q}>1$. Since $\log t$ is concave, Jensen's inequality with respect to the probability measure $d\sigma$ gives
	\begin{equation*}
		\int_X \log  \left[ \left(\frac{g^q}{\int_X g^q\,d\mu}\right)^{\theta-1} \right]d\sigma \le \log \left[ \int_X \left(\frac{g^q}{\int_X g^q\,d\mu}\right)^{\theta-1} d\sigma\right] =\log \left( \frac{\int_X g^p\,d\mu} {\left(\int_X g^q\,d\mu\right)^{p/q}} \right).
	\end{equation*}
	We observe that
	\begin{equation*}
		\int_X \log  \left[ \left(\frac{g^q}{\int_X g^q\,d\mu}\right)^{\theta-1} \right]d\sigma = \frac{p-q}{q}\,\frac{1}{\int_X g^q\,d\mu} \int_X g^q\log \left( \frac{g^q}{\int_X g^q\,d\mu}\right)d\mu.
	\end{equation*}
	Multiplying the previous estimate by $\tfrac{q}{p-q}\bigl(\int_X g^q\,d\mu\bigr)$ and rewriting the logarithm, we obtain the desired result.

	Suppose now that $\bigl(\int_X g^p\,d\mu\bigr)^{\frac{1}{p}}\le A$.  Since $\mu(X)=1$, Jensen's inequality yields
	\begin{equation*}
		0<\int_X g^q\,d\mu \le\left(\int_X g^p\,d\mu\right)^{q/p}\le A^q.
	\end{equation*}
	Finally, the function $t\mapsto t\log(A^q/t)$ on $[0,A^q]$ reaches its maximum at $t=\frac{A^q}{e}$, and hence
	\begin{equation*}
		t\log \left(\frac{A^q}{t}\right)\le\frac{A^q}{e}.
	\end{equation*}
	This proves both inequalities in the second assertion.
\end{proof}

\begin{remark}
	The constant $\frac{p}{p-q}$ in Lemma \ref{lem:probability} cannot be improved in general. Indeed, if $c > 0$, $E \subset X$ is a measurable set with $0 < \mu(E) < 1$, and $g = c\chi_E$, then
	\begin{equation*}
		\int_X g^q\log \left(\frac{g^q}{\int_X g^q\,d\mu}\right)d\mu =c^q\mu(E)\log \left(\frac{1}{\mu(E)}\right),
	\end{equation*}
	whereas
	\begin{equation*}
		\left(\int_X g^q\,d\mu\right) \log \left( \frac{(\int_X g^p\,d\mu)^{q/p}}{\int_X g^q\,d\mu}\right) =c^q\mu(E)\frac{p-q}{p}\log \left(\frac{1}{\mu(E)}\right).
	\end{equation*}
	Thus equality holds. 
\end{remark}

We next recall the self-improving properties of generalized Poincaré inequalities. Let $\mathcal Q$ denote the family of cubes in $\R^n$ with sides parallel to the coordinate axes, and let $a:\mathcal Q\to[0,\infty)$ be a functional. The starting point is an estimate of the form
\begin{equation}\label{eq:generalized-Poincare}
	\frac{1}{|Q|}\int_Q|f(x)-f_Q|\,dx\le a(Q)
\end{equation}
for every cube $Q$, where $f\in L^1_{\loc}(\R^n)$. The question is whether a condition on $a$ allows us to replace the average on the left-hand side by a higher integrability norm of the oscillation.

A general functional approach to this question was developed in \cite{FPW98} and further investigated in \cite{MP98}. This theory studies which geometric summability properties of $a$, over families of pairwise disjoint subcubes of a cube, imply an improvement of the left-hand side of \eqref{eq:generalized-Poincare}. We recall the conditions needed here in the unweighted setting.

\begin{definition}\label{def:Dp}
	Let $1<p<\infty$. We say that $a$ satisfies the $D_p$ condition if there is a constant $C\geq1$ such that
	\begin{equation*}
		\sum_j a(Q_j)^p|Q_j|\le C^p a(Q)^p|Q|
	\end{equation*}
	for every cube $Q$ and every finite family of pairwise disjoint subcubes $\{Q_j\}$ of $Q$. The least such constant is denoted by $\|a\|_{D_p}$, and we write $a\in D_p$.
\end{definition}

Under this condition, \eqref{eq:generalized-Poincare} implies the weak estimate
\begin{equation*}
	\sup_{\lambda>0}\lambda\left(\frac{|\{x\in Q:|f(x)-f_Q|>\lambda\}|}{|Q|}\right)^{\frac{1}{p}}\le C_n\,p\,\|a\|_{D_p}a(Q)
\end{equation*}
where $C_n$ depends only on $n$. This is the unweighted case of \cite{CP}*{Theorem~1.6}, which refines the self-improving result of \cite{FPW98} (see also \cite{PR}*{Theorem~1.2}) and gives a constant linear in $p$. The expression on the left-hand side is the normalized weak $L^p$ quasi-norm of $f-f_Q$ on $Q$.

A stronger summability condition, which also keeps track of the measure occupied by the subcubes, was introduced in \cite{PR}. This smallness condition yields a strong $L^p$ estimate directly.

\begin{definition}\label{def:SD}
	Let $1<p<\infty$ and $1<s<\infty$. We say that $a$ satisfies the $SD_p^s$ condition if there is a constant $C\geq1$ such that
	\begin{equation*}
		\sum_j a(Q_j)^p|Q_j|\le C^p\left(\frac{|\bigcup_jQ_j|}{|Q|}\right)^{\frac{p}{s}}a(Q)^p |Q|
	\end{equation*}
	for every cube $Q$ and every finite family of pairwise disjoint subcubes $\{Q_j\}$ of $Q$. The least such constant is denoted by $\|a\|_{SD_p^s}$, and we write $a\in SD_p^s$.
\end{definition}

Since  $|\bigcup_jQ_j|\le|Q|$, we have $SD_p^s\subset D_p$. The additional factor requires the corresponding sum to become small when the subcubes occupy a small part of $Q$. A basic example, for $n\geq2$, is
\begin{equation*}
	a(Q)=\ell(Q)\left(\frac{1}{|Q|}\int_Q h(x)^p\,dx\right)^{\frac{1}{p}}
\end{equation*}
where $h\geq0$ and $h\in L^p_{\loc}(\R^n)$. This functional belongs to $SD_p^n$ with $\|a\|_{SD_p^n}=1$. The choice $h=|\nabla f|$ gives the usual Poincaré functional.

The condition $SD_p^s$ also allows further gains of integrability. In the weighted setting, weak estimates at Sobolev exponents and exponential integrability in the limiting and supercritical ranges were obtained in \cite{ClarosPerezZheng}. For our logarithmic estimates, it is enough to use strong $L^p$ control of the oscillation and then apply Lemma \ref{lem:probability} with a smaller exponent.

The result in \cite{PR}*{Theorem~1.5} provides a strong estimate with constant $C_n s\|a\|_{SD_p^s}^s$. Lerner, Lorist, and Ombrosi \cite{LLO}*{Theorem~5.3} subsequently obtained a linear dependence on $\|a\|_{SD_p^s}$ by means of sparse domination. In the unweighted case, their estimate reads as follows.

\begin{theorem}\label{thm:strong-selfimproving}
	Let $1<p<\infty$ and $1<s<\infty$, and assume that $a\in SD_p^s$. Let $f\in L^1_{\loc}(\R^n)$ satisfy \eqref{eq:generalized-Poincare} for every cube $Q$. Then
	\begin{equation}\label{eq:strongPS}
		\left(\frac{1}{|Q|}\int_Q|f(x)-f_Q|^p\,dx\right)^{\frac{1}{p}}\le\widetilde a(Q)
	\end{equation}
	for every cube $Q$, where
	\begin{equation*}
		\widetilde a(Q):=C_n(s+1)\|a\|_{SD_p^s}a(Q).
	\end{equation*}
\end{theorem}

Combining Theorem~\ref{thm:strong-selfimproving} with Lemma~\ref{lem:probability}, we obtain a logarithmic form of the self-improving property. The proof only uses the strong estimate \eqref{eq:strongPS}, so $\widetilde a$ may be any nonnegative functional for which it holds.

\begin{theorem}\label{thm:log-selfimproving}
	Let $1<p<\infty$ and let $f\in L^1_{\loc}(\R^n)$. Assume that \eqref{eq:strongPS} holds for every cube $Q$ with a functional $\widetilde a:\mathcal Q\to[0,\infty)$. For every $0<q<p$ and every cube $Q$ on which $f$ is nonconstant, we have
	\begin{align*}
		&\frac{1}{|Q|}\int_Q|f(x)-f_Q|^q\log \left(\frac{|f(x)-f_Q|^q}{\frac{1}{|Q|}\int_Q|f(y)-f_Q|^q\,dy}\right)dx\\
		&\quad\le\frac{p}{p-q}\left(\frac{1}{|Q|}\int_Q|f(x)-f_Q|^q\,dx\right)\log \left(\frac{\widetilde a(Q)^q}{\frac{1}{|Q|}\int_Q|f(y)-f_Q|^q\,dy}\right)\\
		&\quad\le\frac{p}{e(p-q)}\widetilde a(Q)^q.
	\end{align*}
\end{theorem}

\begin{proof}
	Fix a cube $Q$. We apply Lemma \ref{lem:probability} on $X=Q$ with the probability measure $d\mu=\frac{dx}{|Q|}$ and $g(x)=|f(x)-f_Q|$. By \eqref{eq:strongPS}, we have
	\begin{equation*}
		\left(\frac{1}{|Q|}\int_Q|f(x)-f_Q|^p\,dx\right)^{\frac{1}{p}}\le\widetilde a(Q).
	\end{equation*}
	Hence, we obtain the desired estimate from the second part of Lemma \ref{lem:probability} with $A=\widetilde a(Q)$.
\end{proof}

\section{Logarithmic Poincaré--Sobolev inequalities}\label{sec:classical}

In this section, we prove Theorem \ref{thm:intro-local} using the classical Poincaré--Sobolev inequality. We then show that the local estimate implies a global logarithmic Sobolev inequality. The argument is close in spirit to that of Bonforte and Grillo \cite{BG05}*{Proposition~2.1}, who deduced a logarithmic Sobolev inequality from the Sobolev inequality, although they work on compact Riemannian manifolds.

\begin{proof}[Proof of Theorem \ref{thm:intro-local}]
	Let $1\le p<n$, let $Q\subset\R^n$ be a cube, and let $f\in W^{1,p}(Q)$ be nonconstant. Set $p^*=\frac{np}{n-p}$. By the classical Poincaré--Sobolev inequality on cubes, we have
	\begin{equation*}
		\left(\frac{1}{|Q|}\int_Q|f(x)-f_Q|^{p^*}\,dx\right)^{\frac{1}{p^*}}\le C_{n,p}\ell(Q)\left(\frac{1}{|Q|}\int_Q|\nabla f(x)|^p\,dx\right)^{\frac{1}{p}}.
	\end{equation*}
	We apply Lemma \ref{lem:probability} to $g=|f-f_Q|$ on the probability space $(Q,\frac{dx}{|Q|})$, with the exponents $p^*$ and $p$. Since
	\begin{equation*}
		\frac{p^*}{p^*-p}=\frac np,
	\end{equation*}
	the second assertion of the lemma gives \eqref{eq:intro-local}.
\end{proof}

The factor $\ell(Q)^p$ in the local inequality is also the one needed to recover an estimate on the whole space. The following proposition makes this passage precise.

\begin{proposition}[Local implies global]\label{prop:local-implies-global}
	Let $1<p<n$ and let $f\in C_c^\infty(\R^n)$ be nonzero.  Then
	\begin{equation}\label{eq:global}
	\begin{aligned}
		&\int_{\R^n}|f(x)|^p\log \left(\frac{|f(x)|^p}{\int_{\R^n}|f(y)|^p\,dy}\right)dx\le\frac np\left(\int_{\R^n}|f|^p\,dx\right)\log \left(\frac{C_{n,p}^p\int_{\R^n}|\nabla f(x)|^p\,dx}{\int_{\R^n}|f(y)|^p\,dy}\right).
	\end{aligned}
	\end{equation}
\end{proposition}

\begin{proof}
	For $R>0$, let $Q_R=(-\frac{R}{2},\frac{R}{2})^n$ and choose $R_0\geq1$ such that $\supp f\subset Q_R$ for every $R\geq R_0$.  For such $R$, we have
	\begin{equation*}
		f_{Q_R}=\frac{1}{|Q_R|}\int_{Q_R}f(x)\,dx=R^{-n}\int_{\R^n}f(x)\,dx.
	\end{equation*}
	Applying the local estimate \eqref{eq:intro-local} on $Q_R$ and multiplying by $|Q_R|=R^n$, we obtain 
	\begin{align*}
		&\int_{Q_R}|f-f_{Q_R}|^p\log \left(\frac{R^n|f-f_{Q_R}|^p}{\int_{Q_R}|f(y)-f_{Q_R}|^p\,dy}\right)dx\\
		&\quad\le\frac np\left(\int_{Q_R}|f(x)-f_{Q_R}|^p\,dx\right)\log \left(\frac{C_{n,p}^pR^p\int_{\R^n}|\nabla f(x)|^p\,dx}{\int_{Q_R}|f(y)-f_{Q_R}|^p\,dy}\right).
	\end{align*}
	Splitting the powers of $R$ inside the two logarithms yields the same additional term
	\begin{equation*}
		n\left(\int_{Q_R}|f(x)-f_{Q_R}|^p\,dx\right)\log R
	\end{equation*}
	on both sides.  We can cancel it, and we obtain
	\begin{align*}
		&\int_{Q_R}|f-f_{Q_R}|^p\log \left(\frac{|f-f_{Q_R}|^p}{\int_{Q_R}|f(y)-f_{Q_R}|^p\,dy}\right)dx\\
		&\quad\le\frac np\left(\int_{Q_R}|f(x)-f_{Q_R}|^p\,dx\right)\log \left(\frac{C_{n,p}^p\int_{\R^n}|\nabla f(x)|^p\,dx}{\int_{Q_R}|f(y)-f_{Q_R}|^p\,dy}\right).
	\end{align*}

	We now pass to the limit.  Since $f_{Q_R}=R^{-n}\int_{\R^n}f$, we have $f_{Q_R}\to 0$, and $(f-f_{Q_R})\chi_{Q_R}\to f$ pointwise on $\R^n$.  By Fatou's lemma,
	\begin{equation*}
		\int_{\R^n}|f(x)|^p\,dx\le\liminf_{R\to\infty}\int_{Q_R}|f(x)-f_{Q_R}|^p\,dx.
	\end{equation*}
	On the other hand, by the triangle inequality and the inclusion $\supp f\subset Q_R$, we have
	\begin{align*}
		\left(\int_{Q_R}|f(x)-f_{Q_R}|^p\,dx\right)^{\frac{1}{p}}&\le\left(\int_{\R^n}|f(x)|^p\,dx\right)^{\frac{1}{p}}+|f_{Q_R}||Q_R|^{\frac{1}{p}}\\
		&=\left(\int_{\R^n}|f(x)|^p\,dx\right)^{\frac{1}{p}}+\left|\int_{\R^n}f(x)\,dx\right|R^{-n(1-\frac{1}{p})}.
	\end{align*}
	Taking the $\limsup$ and combining the resulting estimate with Fatou's lemma, we obtain
	\begin{equation}\label{eq:limQR}
		\int_{Q_R}|f(x)-f_{Q_R}|^p\,dx\longrightarrow\int_{\R^n}|f(x)|^p\,dx>0.
	\end{equation}

	It remains to pass to the limit in the logarithmic part. We first consider the integral over the fixed set $\supp f$.  For every $x\in\supp f$, we have $f(x)-f_{Q_R} \to f(x)$, and the denominator in the logarithm satisfies \eqref{eq:limQR}. Hence, the integrand converges pointwise on $\supp f$ to the desired limiting integrand.  To justify the application of the dominated convergence theorem, we write
	\begin{align*}
		&|f(x)-f_{Q_R}|^p\log \left(\frac{|f(x)-f_{Q_R}|^p}{\int_{Q_R}|f(y)-f_{Q_R}|^p\,dy}\right)\\
		&\quad=|f(x)-f_{Q_R}|^p\log|f(x)-f_{Q_R}|^p-|f(x)-f_{Q_R}|^p\log \left(\int_{Q_R}|f(y)-f_{Q_R}|^p\,dy\right).
	\end{align*}
	For all sufficiently large $R$, we have $|f_{Q_R}|\le1$, and hence $|f(x)-f_{Q_R}|\le\|f\|_{L^\infty(\R^n)}+1$, while \eqref{eq:limQR} implies that the logarithm of the denominator is uniformly bounded.  Since $t\mapsto|t|^p\log|t|^p$ is continuous and therefore bounded on every bounded interval, the integrand is bounded by a constant independent of $R$ and $x\in\supp f$.  Since $\supp f$ has finite measure, this constant is an integrable dominating function.  The dominated convergence theorem now gives
	\begin{align*}
		&\int_{\supp f}|f(x)-f_{Q_R}|^p\log \left(\frac{|f(x)-f_{Q_R}|^p}{\int_{Q_R}|f(y)-f_{Q_R}|^p\,dy}\right)dx\\
		&\qquad\longrightarrow\int_{\R^n}|f(x)|^p\log \left(\frac{|f(x)|^p}{\int_{\R^n}|f(y)|^p\,dy}\right)dx.
	\end{align*}
	We next consider $Q_R\setminus\supp f$.  Since $f=0$ there, the integrand is constant as a function of $x$.  If $\int_{\R^n}f=0$, then $f_{Q_R}=0$, and this contribution is zero.  Suppose that $\int_{\R^n}f\neq0$.  By \eqref{eq:limQR}, for each sufficiently large $R$, we have
	\begin{equation*}
		\frac{1}{2}\int_{\R^n}|f(y)|^p\,dy\le\int_{Q_R}|f(y)-f_{Q_R}|^p\,dy\le2\int_{\R^n}|f(y)|^p\,dy.
	\end{equation*}
	Since $f_{Q_R}=R^{-n}\int_{\R^n}f$, it follows that
	\begin{align*}
		\left|\log \left(\frac{|f_{Q_R}|^p}{\int_{Q_R}|f(y)-f_{Q_R}|^p\,dy}\right)\right|&=\left|\log \left(\frac{|\int_{\R^n}f(x)\,dx|^p}{\int_{Q_R}|f(y)-f_{Q_R}|^p\,dy}\right)-np\log R\right|\\
		&\le \left|\log \left(\frac{|\int_{\R^n}f(x)\,dx|^p}{\int_{\R^n}|f(y)|^p\,dy}\right)\right|+\log2+np\log R\\
		&=: C+np\log R.
	\end{align*}
	Consequently,
	\begin{align*}
		&\left|\int_{Q_R\setminus\supp f}|f_{Q_R}|^p\log \left(\frac{|f_{Q_R}|^p}{\int_{Q_R}|f(y)-f_{Q_R}|^p\,dy}\right)dx\right|\le|Q_R||f_{Q_R}|^p(C+np\log R)\\
		&\quad=\left|\int_{\R^n}f(x)\,dx\right|^pR^{-n(p-1)}(C+np\log R)\longrightarrow0,
	\end{align*}
	because $p>1$. Combining the two limits just proved, we conclude that
	$$\int_{Q_R}|f(x)-f_{Q_R}|^p\log \left(\frac{|f(x)-f_{Q_R}|^p}{\int_{Q_R}|f(y)-f_{Q_R}|^p\,dy}\right)dx \longrightarrow\int_{\R^n}|f(x)|^p\log \left(\frac{|f(x)|^p}{\int_{\R^n}|f(y)|^p\,dy}\right)dx.$$
	It remains to pass to the limit on the right-hand side. Since $f$ is nonzero and compactly supported, we have $\|\nabla f\|_{L^p(\R^n)}>0$. Hence the function
	\begin{equation*}
		t\longmapsto\frac np t\log \left(\frac{C_{n,p}^p\int_{\R^n}|\nabla f(x)|^p\,dx}{t}\right)
	\end{equation*}
	is continuous for $t>0$.  Since the limit \eqref{eq:limQR} is strictly positive, this continuity gives
	\begin{align*}
		&\frac np\left(\int_{Q_R}|f(x)-f_{Q_R}|^p\,dx\right)\log \left(\frac{C_{n,p}^p\int_{\R^n}|\nabla f(x)|^p\,dx}{\int_{Q_R}|f(y)-f_{Q_R}|^p\,dy}\right)\\
		&\qquad\longrightarrow\frac np\left(\int_{\R^n}|f(x)|^p\,dx\right)\log \left(\frac{C_{n,p}^p\int_{\R^n}|\nabla f(x)|^p\,dx}{\int_{\R^n}|f(y)|^p\,dy}\right).
	\end{align*}
	We may therefore pass to the limit in the inequality, which yields \eqref{eq:global} and completes the proof.
\end{proof}

The restriction $p>1$ is used in the estimates for the contribution of the average $f_{Q_R}$ outside the support of $f$. In particular, it ensures that $R^{-n(p-1)}(1+\log R)$ tends to zero. Thus, this argument does not give the global conclusion at $p=1$, although the local inequality in Theorem \ref{thm:intro-local} includes that endpoint.

Since the previous inequality is homogeneous with respect to multiplication of $f$ by scalars, we may impose the normalization
\begin{equation*}
	\int_{\mathbb R^n}|f(x)|^p\,dx=1.
\end{equation*}
Under this normalization, Proposition \ref{prop:local-implies-global} becomes
\begin{equation*}
	\int_{\mathbb R^n}|f(x)|^p\log|f(x)|^p\,dx\le\frac np\log \left(C_{n,p}^p\int_{\mathbb R^n}|\nabla f(x)|^p\,dx\right).
\end{equation*}
This is the Euclidean $L^p$ logarithmic Sobolev inequality \eqref{eq:log-global-Lp} with $C_{n,p}^p$ in place of the optimal constant $\alpha_{n,p}$, and with the coefficient $\frac np$ given by scaling (Remark~\ref{rem:scaling}).

\section{Fractional logarithmic Poincaré--Sobolev inequalities}\label{sec:fractional}

In this section, we apply the same argument to the fractional Poincaré--Sobolev inequality. After proving Theorem \ref{thm:intro-fractional}, we pass from the local estimate to the whole space. The factor $1-\delta$ gives the connection with the first-order result as $\delta\to1^-$.

\begin{proof}[Proof of Theorem \ref{thm:intro-fractional}]
	Let $0<\delta<1$, let $1\le p<\frac{n}{\delta}$, and set $p_\delta^*=\frac{np}{n-\delta p}$. Fix a cube $Q$ and let $f\in W^{\delta,p}(Q)$ be nonconstant. The fractional Poincaré--Sobolev inequality in \cite{MPW}*{Theorem~5.7}, with Lebesgue measure, gives
	\begin{align*}
		&\left(\frac{1}{|Q|}\int_Q|f(x)-f_Q|^{p_\delta^*}\,dx\right)^{\frac{1}{p_\delta^*}}\le C_{n,p,\delta}(1-\delta)^{\frac{1}{p}}\ell(Q)^\delta\left(\frac{1}{|Q|}\int_Q\int_Q\frac{|f(x)-f(y)|^p}{|x-y|^{n+\delta p}}\,dy\,dx\right)^{\frac{1}{p}}.
	\end{align*}
	We now apply Lemma \ref{lem:probability} to $g=|f-f_Q|$ on $(Q,\frac{dx}{|Q|})$, with the exponents $p_\delta^*$ and $p$. Since
	\begin{equation*}
		\frac{p_\delta^*}{p_\delta^*-p}=\frac{n}{\delta p},
	\end{equation*}
	we obtain \eqref{eq:intro-fractional}.
\end{proof}

We next pass from the local fractional estimate to an inequality on the whole space.

\begin{proposition}[Fractional local implies global]\label{prop:fractional-local-implies-global}
	Let $0<\delta<1$, let $1<p<\frac{n}{\delta}$, and let $f\in C_c^\infty(\R^n)$ be nonzero. Then
	\begin{equation}\label{eq:global-fractional}
	\begin{aligned}
		&\int_{\R^n}|f(x)|^p\log \left(\frac{|f(x)|^p}{\int_{\R^n}|f(y)|^p\,dy}\right)dx\\
		&\quad\le\frac{n}{\delta p}\left(\int_{\R^n}|f(x)|^p\,dx\right)\log \left(\frac{C_{n,p,\delta}^p(1-\delta)\int_{\R^n}\int_{\R^n}\frac{|f(x)-f(y)|^p}{|x-y|^{n+\delta p}}\,dy\,dx}{\int_{\R^n}|f(y)|^p\,dy}\right).
	\end{aligned}
	\end{equation}
\end{proposition}

\begin{proof}
	For $R>0$, let $Q_R=(-\frac{R}{2},\frac{R}{2})^n$ and choose $R_0\geq1$ such that $\supp f\subset Q_R$ for every $R\geq R_0$. For such $R$, we have
	\begin{equation*}
		f_{Q_R}=\frac{1}{|Q_R|}\int_{Q_R}f(x)\,dx=R^{-n}\int_{\R^n}f(x)\,dx.
	\end{equation*}
	Applying the local estimate \eqref{eq:intro-fractional} on $Q_R$ and multiplying by $|Q_R|=R^n$, we obtain
	\begin{align*}
		&\int_{Q_R}|f-f_{Q_R}|^p\log \left(\frac{R^n|f-f_{Q_R}|^p}{\int_{Q_R}|f(y)-f_{Q_R}|^p\,dy}\right)dx\\
		&\quad\le\frac{n}{\delta p}\left(\int_{Q_R}|f(x)-f_{Q_R}|^p\,dx\right)\log \left(\frac{C_{n,p,\delta}^p(1-\delta)R^{\delta p}\int_{Q_R}\int_{Q_R}\frac{|f(x)-f(y)|^p}{|x-y|^{n+\delta p}}\,dy\,dx}{\int_{Q_R}|f(y)-f_{Q_R}|^p\,dy}\right).
	\end{align*}
	Splitting the power of $R$ inside the logarithm on the left-hand side produces the additional term
	\begin{equation*}
		n\left(\int_{Q_R}|f(x)-f_{Q_R}|^p\,dx\right)\log R.
	\end{equation*}
	On the right-hand side, we have $\frac{n}{\delta p}\log R^{\delta p}=n\log R$, so the same additional term appears. We can cancel it, and we obtain
	\begin{align*}
		&\int_{Q_R}|f-f_{Q_R}|^p\log \left(\frac{|f-f_{Q_R}|^p}{\int_{Q_R}|f(y)-f_{Q_R}|^p\,dy}\right)dx\\
		&\quad\le\frac{n}{\delta p}\left(\int_{Q_R}|f(x)-f_{Q_R}|^p\,dx\right) \log \left(\frac{C_{n,p,\delta}^p(1-\delta)\int_{Q_R}\int_{Q_R}\frac{|f(x)-f(y)|^p}{|x-y|^{n+\delta p}}\,dy\,dx}{\int_{Q_R}|f(y)-f_{Q_R}|^p\,dy}\right).
	\end{align*}

	The convergence arguments for the $L^p$ norms and for the logarithmic integrals in the proof of Proposition \ref{prop:local-implies-global} only use $p>1$ and $f\in C_c^\infty(\R^n)$, and do not require $p<n$. They therefore apply throughout the range $1<p<\frac{n}{\delta}$. In particular, as in \eqref{eq:limQR},
	\begin{equation*}
		\int_{Q_R}|f(x)-f_{Q_R}|^p\,dx\longrightarrow\int_{\R^n}|f(x)|^p\,dx>0.
	\end{equation*}
	The argument for the logarithmic integral is also unchanged. On $\supp f$, the dominated convergence theorem applies because $f_{Q_R}\to0$, the functions $f-f_{Q_R}$ are uniformly bounded, and the denominator converges to a strictly positive number. On $Q_R\setminus\supp f$, the absolute value of the corresponding integral is bounded by
	\begin{equation*}
		CR^{-n(p-1)}(1+\log R),
	\end{equation*}
	which tends to zero because $p>1$. Consequently,
	\begin{align*}
		&\int_{Q_R}|f(x)-f_{Q_R}|^p\log \left(\frac{|f(x)-f_{Q_R}|^p}{\int_{Q_R}|f(y)-f_{Q_R}|^p\,dy}\right)dx \longrightarrow\int_{\R^n}|f(x)|^p\log \left(\frac{|f(x)|^p}{\int_{\R^n}|f(y)|^p\,dy}\right)dx.
	\end{align*}

	It remains to pass to the limit in the fractional seminorm. Since the cubes $Q_R$ increase to $\R^n$ and the integrand is nonnegative, the monotone convergence theorem gives
	\begin{equation}\label{eq:fractional-limit}
		\int_{Q_R}\int_{Q_R}\frac{|f(x)-f(y)|^p}{|x-y|^{n+\delta p}}\,dy\,dx\longrightarrow\int_{\R^n}\int_{\R^n}\frac{|f(x)-f(y)|^p}{|x-y|^{n+\delta p}}\,dy\,dx.
	\end{equation}
	The fractional seminorm on the right-hand side is finite. Indeed, $f\in C_c^\infty(\R^n)\subset W^{1,p}(\R^n)$, and $W^{1,p}(\R^n)\subset W^{\delta,p}(\R^n)$ continuously (see \cite{Hitchhiker}*{Proposition~2.2 and its proof}). The double integral in \eqref{eq:fractional-limit} is strictly positive. Otherwise, $f(x)=f(y)$ for almost every pair $(x,y)\in\R^n\times\R^n$, so $f$ would be constant almost everywhere. Since $f$ is compactly supported, this implies that $f=0$, contrary to the hypothesis. Combining \eqref{eq:limQR}, \eqref{eq:fractional-limit}, and the continuity of the logarithm gives
	\begin{align*}
		&\frac{n}{\delta p}\left(\int_{Q_R}|f(x)-f_{Q_R}|^p\,dx\right)\log \left(\frac{C_{n,p,\delta}^p(1-\delta)\int_{Q_R}\int_{Q_R}\frac{|f(x)-f(y)|^p}{|x-y|^{n+\delta p}}\,dy\,dx}{\int_{Q_R}|f(y)-f_{Q_R}|^p\,dy}\right)\\
		&\qquad\longrightarrow\frac{n}{\delta p}\left(\int_{\R^n}|f(x)|^p\,dx\right)\log \left(\frac{C_{n,p,\delta}^p(1-\delta)\int_{\R^n}\int_{\R^n}\frac{|f(x)-f(y)|^p}{|x-y|^{n+\delta p}}\,dy\,dx}{\int_{\R^n}|f(y)|^p\,dy}\right).
	\end{align*}
	We may therefore pass to the limit in the inequality, which yields \eqref{eq:global-fractional} and completes the proof.
\end{proof}

Since the previous inequality is homogeneous with respect to multiplication of $f$ by scalars, we may impose the normalization
\begin{equation*}
	\int_{\R^n}|f(x)|^p\,dx=1.
\end{equation*}
Under this normalization, Proposition \ref{prop:fractional-local-implies-global} becomes
\begin{equation*}
	\int_{\R^n}|f(x)|^p\log|f(x)|^p\,dx\le\frac{n}{\delta p}\log \left(C_{n,p,\delta}^p(1-\delta)\int_{\R^n}\int_{\R^n}\frac{|f(x)-f(y)|^p}{|x-y|^{n+\delta p}}\,dy\,dx\right).
\end{equation*}
As explained in Remark~\ref{rem:fractional-scaling}, the coefficient $\frac{n}{\delta p}$ is the one given by scaling. The factor $1-\delta$ is the Bourgain--Brezis--Mironescu normalization \cite{BBM2}, which connects the fractional seminorm with the $L^p$ norm of the gradient $|\nabla f|$ as $\delta\to1^-$. We do not claim that the resulting constant is optimal.

\section{Local Nash inequalities}\label{sec:nash}

Following the argument of Beckner \cite{Beckner}, we combine the local logarithmic inequality with another application of Jensen's inequality to obtain a local $L^p$ Nash inequality.

\begin{proposition}\label{prop:local-Nash}
	Let $1<p<n$, let $Q\subset\R^n$ be a cube, and let $f\in W^{1,p}(Q)$. Then
	\begin{align*}
		&\left(\frac{1}{|Q|}\int_Q|f-f_Q|^p\,dx\right)^{1+\frac{p}{n(p-1)}}\le C_{n,p}^p\ell(Q)^p\left(\frac{1}{|Q|}\int_Q|\nabla f|^p\,dx\right)\left(\frac{1}{|Q|}\int_Q|f-f_Q|\,dx\right)^{\frac{p^2}{n(p-1)}}.
	\end{align*}
\end{proposition}

\begin{proof}
If $f = f_Q$ almost everywhere on $Q$, the inequality holds trivially. Otherwise, by homogeneity, it suffices to assume that
	\begin{equation*}
		\frac{1}{|Q|}\int_Q|f(x)-f_Q|^p\,dx=1.
	\end{equation*}
	Applying \eqref{eq:intro-local} gives
	\begin{equation*}
		\frac{1}{|Q|}\int_Q|f(x)-f_Q|^p\log \left(|f(x)-f_Q|^p\right)dx\le\frac np\log \left(C_{n,p}^p\ell(Q)^p\frac{1}{|Q|}\int_Q|\nabla f(x)|^p\,dx\right).
	\end{equation*}
	Dividing by $p$, we obtain
	\begin{equation*}
		\frac{1}{|Q|}\int_Q|f(x)-f_Q|^p\log|f(x)-f_Q|\,dx\le\frac{n}{p^2}\log \left(C_{n,p}^p\ell(Q)^p\frac{1}{|Q|}\int_Q|\nabla f(x)|^p\,dx\right).
	\end{equation*}
	The measure
	\begin{equation*}
		d\nu(x)=\frac{|f(x)-f_Q|^p}{|Q|}\,dx
	\end{equation*}
	is a probability measure on $Q$. Let $E=\{x\in Q:|f(x)-f_Q|>0\}$. Then $\nu(E)=1$, and the function $1/|f-f_Q|^{p-1}$ is positive and integrable with respect to $\nu$ on $E$. The logarithmic integral is finite by Lemma \ref{lem:probability} and the classical Poincaré--Sobolev inequality. Applying Jensen's inequality to the concave function $\log t$ gives
	\begin{align*}
		-(p-1)\frac{1}{|Q|}\int_Q|f(x)-f_Q|^p\log|f(x)-f_Q|\,dx& =\int_E\log \left(\frac{1}{|f(x)-f_Q|^{p-1}}\right)d\nu(x)\\
		&\le\log \left(\int_E\frac{1}{|f(x)-f_Q|^{p-1}}\,d\nu(x)\right)\\
		&=\log \left(\frac{1}{|Q|}\int_E\frac{|f(x)-f_Q|^p}{|f(x)-f_Q|^{p-1}}\,dx\right)\\
		&=\log \left(\frac{1}{|Q|}\int_Q|f(x)-f_Q|\,dx\right).
	\end{align*}
	Equivalently,
	\begin{equation*}
		-\frac{1}{p-1}\log \left(\frac{1}{|Q|}\int_Q|f(x)-f_Q|\,dx\right)\le\frac{1}{|Q|}\int_Q|f(x)-f_Q|^p\log|f(x)-f_Q|\,dx.
	\end{equation*}
	Combining the last two estimates gives
	\begin{equation*}
		-\frac{p^2}{n(p-1)}\log \left(\frac{1}{|Q|}\int_Q|f(x)-f_Q|\,dx\right)\le\log \left(C_{n,p}^p\ell(Q)^p\frac{1}{|Q|}\int_Q|\nabla f(x)|^p\,dx\right).
	\end{equation*}
	Exponentiating, we obtain
	\begin{equation*}
		1\le C_{n,p}^p\ell(Q)^p\left(\frac{1}{|Q|}\int_Q|\nabla f(x)|^p\,dx\right)\left(\frac{1}{|Q|}\int_Q|f(x)-f_Q|\,dx\right)^{\frac{p^2}{n(p-1)}}.
	\end{equation*}
	Removing the normalization establishes the proposition.
\end{proof}

Taking $p=2$ in Proposition~\ref{prop:local-Nash}, we immediately obtain the following local version of the classical Nash inequality.

\begin{corollary}[Local Nash inequality]\label{cor:local-Nash-L2}
	Let $n\geq3$, let $Q\subset\R^n$ be a cube, and let $f\in W^{1,2}(Q)$. Then
	\begin{equation*}
		\left(\frac{1}{|Q|}\int_Q|f-f_Q|^2\,dx\right)^{1+2/n}\le C_{n,2}^2\ell(Q)^2\left(\frac{1}{|Q|}\int_Q|\nabla f|^2\,dx\right)\left(\frac{1}{|Q|}\int_Q|f-f_Q|\,dx\right)^{4/n}.
	\end{equation*}
\end{corollary}

The same argument applied to Theorem~\ref{thm:intro-fractional} gives a fractional version. Global fractional Gagliardo--Nirenberg inequalities were recently studied in \cite{HLYY26}.

\begin{proposition}\label{prop:fractional-Nash}
	Let $0<\delta<1$, let $1<p<n/\delta$, let $Q\subset\R^n$ be a cube, and let $f\in W^{\delta,p}(Q)$. Then
	\begin{align*}
		&\left(\frac{1}{|Q|}\int_Q|f-f_Q|^p\,dx\right)^{1+\frac{\delta p}{n(p-1)}}\\
		&\quad\le C_{n,p,\delta}^p(1-\delta)\ell(Q)^{\delta p}\left(\frac{1}{|Q|}\int_Q\int_Q\frac{|f(x)-f(y)|^p}{|x-y|^{n+\delta p}}\,dy\,dx\right)\left(\frac{1}{|Q|}\int_Q|f-f_Q|\,dx\right)^{\frac{\delta p^2}{n(p-1)}}.
	\end{align*}
\end{proposition}

\begin{proof}
	 If $f=f_Q$ almost everywhere on $Q$, the inequality holds trivially. Otherwise, by homogeneity, it suffices to assume that
	\begin{equation*}
		\frac{1}{|Q|}\int_Q|f(x)-f_Q|^p\,dx=1.
	\end{equation*}
	Dividing the estimate in Theorem \ref{thm:intro-fractional} by $p$, we obtain
	\begin{align*}
		&\frac{1}{|Q|}\int_Q|f(x)-f_Q|^p\log|f(x)-f_Q|\,dx\\
		&\quad\le\frac{n}{\delta p^2}\log \left(C_{n,p,\delta}^p(1-\delta)\ell(Q)^{\delta p}\frac{1}{|Q|}\int_Q\int_Q\frac{|f(x)-f(y)|^p}{|x-y|^{n+\delta p}}\,dy\,dx\right).
	\end{align*}
	The argument with the probability measure $d\nu=|f-f_Q|^p\frac{dx}{|Q|}$ in the proof of Proposition \ref{prop:local-Nash} also gives
	\begin{equation*}
		-\frac{1}{p-1}\log \left(\frac{1}{|Q|}\int_Q|f(x)-f_Q|\,dx\right)\le\frac{1}{|Q|}\int_Q|f(x)-f_Q|^p\log|f(x)-f_Q|\,dx.
	\end{equation*}
	Here the finiteness of the logarithmic integral follows from the fractional Poincaré--Sobolev inequality and Lemma \ref{lem:probability}. Combining the two estimates and exponentiating, we obtain
	$$1\le C_{n,p,\delta}^p(1-\delta)\ell(Q)^{\delta p}\left(\frac{1}{|Q|}\int_Q\int_Q\frac{|f(x)-f(y)|^p}{|x-y|^{n+\delta p}}\,dy\,dx\right) \left(\frac{1}{|Q|}\int_Q|f(x)-f_Q|\,dx\right)^{\frac{\delta p^2}{n(p-1)}}.$$
	Finally, applying this normalized estimate to
	\begin{equation*}
		\frac{f}{\left(\frac{1}{|Q|}\int_Q|f(x)-f_Q|^p\,dx\right)^{\frac{1}{p}}}
	\end{equation*}
	yields the desired inequality.
\end{proof}

\subsection*{Acknowledgment}
The authors are very grateful to Prof. Neal Bez for drawing their attention to some references on logarithmic Sobolev inequalities.

\subsection*{Use of AI tools}
Claude Opus 5.5 (Anthropic) was used for language editing, drafting assistance, and checking references and computations. The authors have reviewed the final text and take full responsibility for its content. All mathematical results are due to the authors.

\end{document}